\documentclass[10pt]{article}
\usepackage{latexsym}
\usepackage{wrapfig}
\usepackage{graphicx}
\usepackage{color}
\usepackage{amsmath}
\usepackage{amsfonts}
\usepackage{amssymb}
\usepackage{url}

\definecolor{darkblue}{RGB}{0, 0, 102}
\definecolor{IKB}{RGB}{0, 47, 167}

\usepackage[pdftex,bookmarks,bookmarksopen=true,pdfstartview={FitH},bookmarksnumbered=true,colorlinks=true,citecolor=darkblue,linkcolor=IKB,plainpages=false,pdfpagelabels]{hyperref}

\newtheorem{theorem}{Theorem}
\newtheorem{lemma}[theorem]{Lemma}

\newcommand{\E}{\begin{equation}}
\newcommand{\EE}{\end{equation}}
\newcommand{\QED}{\ \rule{.1in}{.1in}}

\advance \topmargin by -\headheight
\advance \topmargin by -\headsep
\begin{document}

\hspace{.1in}\\

\vskip 1cm

\begin{center}
\Large {\bf Strong NP-hardness of AC power flows feasibility}

\normalsize

Daniel Bienstock and Abhinav Verma, Columbia University
\end{center}
\section{Introduction}
The AC-OPF problem is a nonlinear, nonconvex problem arising in the study
of electrical networks. This problem has lately received increased
attention, spurred by the work in \cite{lavaeilow}.  A recent survey is given in
\cite{danianreview}. An important question
concerns the fundamental complexity of AC-OPF.  In this
setting we note \cite{pascal}, which discusses a
proof of weak NP-hardness of AC-OPF on trees\footnote{Note: the reduction encodes some irrational quantities.}; also see \cite{lavaeilow} for an outline of a proof.  On graphs with bounded tree-width the problem can be solved to any given
tolerance $\epsilon$ in time polynomial in the size of the network and $\epsilon^{-1}$ \cite{gonzaloandi}.

The purpose of this note is to present a
rigorous proof of strong NP-hardness of AC-OPF on general graphs; this proof
builds on a section in the PhD thesis \cite{abhinav}.  A challenge in
the development of such a proof is the fact that the solution to 
an AC-OPF problem may have irrational coordinates; this necessitates an
elaboration in our technique.

The problem we consider is described by a directed graph $G$ representing a power transmission system, where
for each line (i.e., arc) $(i,j)$ we are given a positive parameter $x_{ij}$, the {\em reactance} and a nonnegative value $\theta^{\max}_{ij}$, the \textit{maximum phase angle difference}.   Additionally for each bus (i.e., vertex) $i$ we have a value
$b_i$ indicating the net generation at $i$.  We assume $\sum_i b_i = 0$.
Denoting by $N$ the node-arc incidence matrix of $G$, the lossless AC power flow feasibility
system on $G$ can be written as 
\begin{subequations} \label{lossless}
\begin{eqnarray}
&& N f \ = \ b, \\
&&  \sin(\theta_i - \theta_j) \ = \ x_{ij} f_{ij}, \quad \forall (i,j), \label{ohm}\\
&&  |\theta_i - \theta_j| \ \le \ \theta^{\max}_{ij}, \quad \forall (i,j).
\end{eqnarray}
\end{subequations}
Typically, $\theta^{\max}_{ij} < \pi/2$ for each line $(i,j)$.  System (\ref{lossless}) is a special case of the general AC power flow problem; given a line
$(i,j)$ variable $f_{ij}$ models the real (active) power flowing from $i$ to $j$, 
which could be negative,
and for any bus $i$ variable $\theta_i$ is the phase angle at $i$.  See \cite{andersson1}, \cite{bergenvittal}, \cite{overbyebook}.  This model, and
variants thereof, was considered in \cite{lesieutremezaetal}, \cite{12BV},
\cite{bbc}.  From a power engineering perspective the model assumes zero
resistances and unconstrained reactive power flows and injections. \\

System (\ref{lossless}) is both nonlinear and nonconvex.  In this section we
prove that testing feasibility of such a system is a strongly NP-hard problem,
that is to say it remains NP-hard even if the number of bits in the input data
is polynomially bounded as a function of the number of buses. 

\subsection{Main construction}
\begin{wrapfigure}[10]{r}{0.32\textwidth}
  \begin{center}
    \includegraphics[width=0.28\textwidth, angle=270]{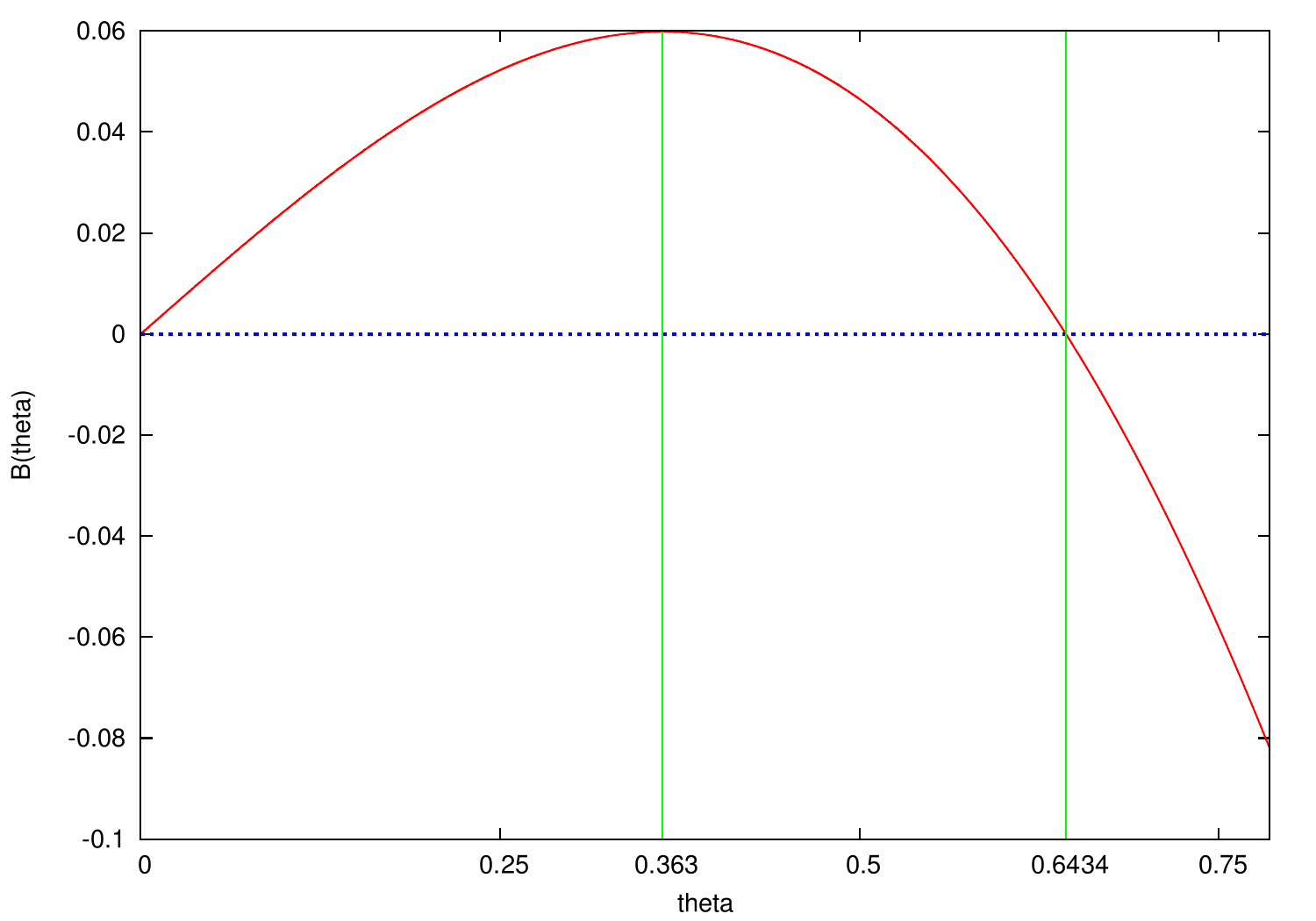}
  \end{center}
\vskip -15pt
  \caption{$\Delta(\theta)$ for $\theta \in [0, \pi/4]$} \label{Bfunction}
\vskip -25pt
\end{wrapfigure}
Define 
$$\Delta(\theta) \ \doteq \ -\sin(\theta) + \frac{5}{8}\sin(2 \theta),$$
and set
$$ \theta_0 \ = \ \cos^{-1}\left(\frac{4}{5}\right) \approx .6435, \quad \theta_1 \ = \ \cos^{-1}\left( \frac{1}{5} + \sqrt{ \frac{1}{25} + \frac{1}{2}}\right) \approx .3630.$$

\noindent Then we have:
\begin{lemma}\label{curvy} Suppose $0 \le \theta \le \frac{\pi}{2}$. Then:
\begin{itemize}
\item [{\bf(a)}] $\Delta(\theta) = 0$ iff $\theta = 0$ or $\theta = \theta_0$.
\item [{\bf(b)}] $\Delta'(\theta) > 0$ for $0 \le \theta < \theta_1$, $\Delta'(\theta_1) = 0$  and $\Delta'(\theta) < 0$ for $\theta_1 < \theta \le \pi/2$.
\end{itemize}
\end{lemma}
\noindent {\em Proof.} Part (a) is clear. To prove (b) let $c = \cos(\theta)$. Then
$\Delta'(\theta) = - c + \frac{5}{4}(2 c^2 - 1)$ whose only zero in $[0,1]$ is at
$1/5 + \sqrt{1/25 + 1/2}$. \QED\\

\noindent See Figure \ref{Bfunction}.

\begin{figure}[ht] 
\centering
\includegraphics[width=0.6\textwidth]{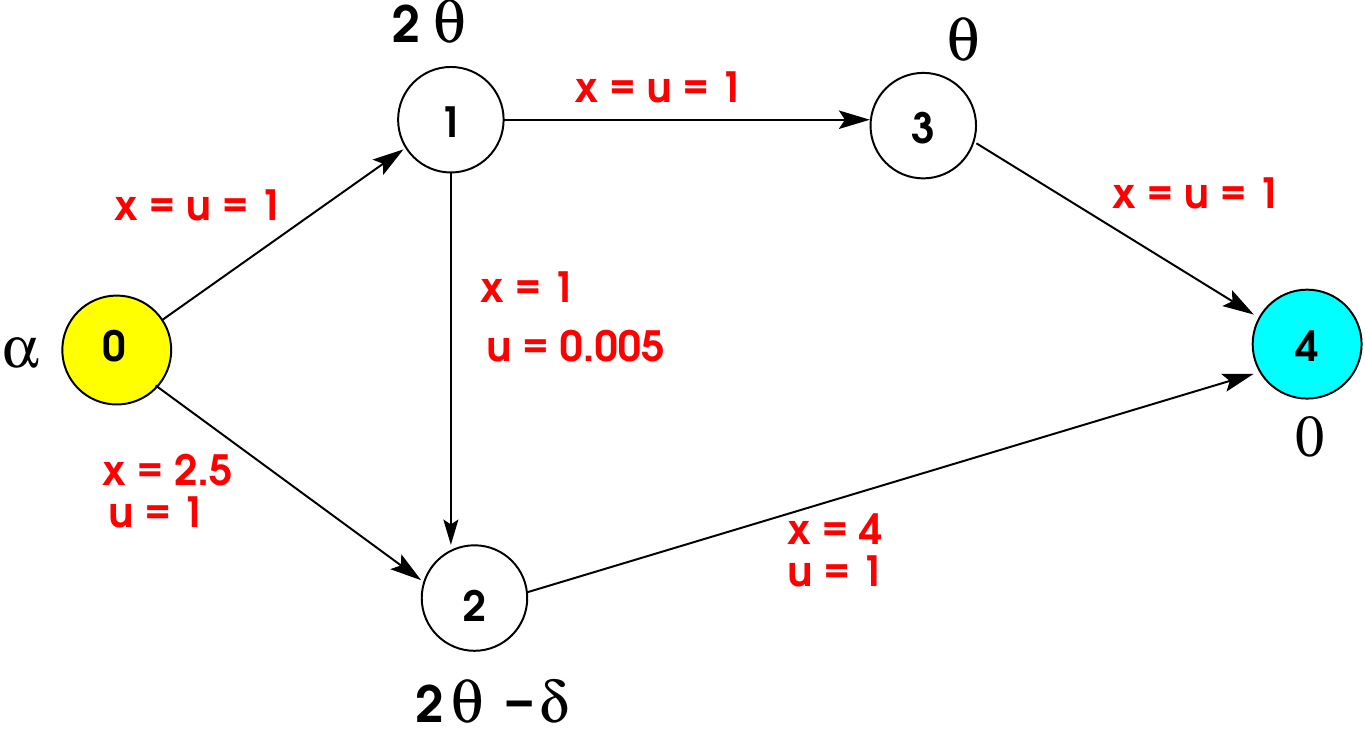}
\vskip -10pt
\caption{Basic transmission system B used in proof} \label{banana}
\end{figure}

Our construction is centered on the transmission system B shown in Figure \ref{banana}.  In this figure we show, next 
to each line, its reactance (``x'') and its limit (``u'').  We consider
solutions to the lossless AC model where bus 0 injects power into the system,
bus 4 withdraws power, and all other buses have zero net balance.  The figure also shows our
naming convention for phase angles in a solution to the lossless AC power flow problem on this graph.  Denoting the phase angle at bus $i$ by $\theta_i$, for $0 \le i \le 4$, we have:

\begin{itemize}
\item[(i)] Without loss of generality $\theta_4 = 0$,
\item[(ii)] We simplify $\theta_3$ as $\theta$.  
\item [(iii)] $\theta_1 = 2 \theta$ (flow conservation at bus $3$).
\item [(iv)] We write $\theta_2$ as $2\theta - \delta$.  Possibly $\delta < 0$.
\item [(v)] We write $\theta_0 = \alpha$.  

\end{itemize}
We can make some basic observations:\\
\noindent {\bf 1.} Since power flows from $3$ to $4$, 
\begin{eqnarray}
&& 0 \le \theta \le \pi/2. \label{basic1}
\end{eqnarray}
\noindent {\bf 2.} The absolute value of the flow on line $(1,2)$ is $|\sin(\delta)|$.
But the flow limit on line $(1,2)$ is $0.005$, implying
\begin{eqnarray}
&& | \sin \delta | \le 0.005, \quad \mbox{and consequently} \quad |\delta| < 0.0050001. \label{basic2}
\end{eqnarray}
\noindent {\bf Remark}.  The bounds in \eqref{basic2} can be made arbitrarily
small by choosing a small enough limit on line (1,2).  \\

\noindent {\bf 3.} Flow on $(0,1)$ and $(0,2)$ must be nonnegative.  Together
with the phase angle limits we have 
\begin{eqnarray}
&& \max\{0, -\delta\}  \le \alpha - 2 \theta \le \min \{\pi/2, \pi/2 - \delta\} \label{basic3}
\end{eqnarray}
\noindent {\bf 4.} Applying a similar reasoning to line $(2,4)$ we get
\begin{eqnarray}
&& \max\{0, \delta\}  \le 2 \theta \le \min \{\pi/2, \pi/2 + \delta\} \label{basic4}
\end{eqnarray}
\noindent {\bf 5.} The flow conservation equations at buses $1$ and $2$ are,
respectively,
\begin{subequations}\label{cons12}
\begin{eqnarray}
\sin(\alpha - 2 \theta) & = & \sin \delta + \sin \theta, \label{cons1}\\
\frac{1}{2.5} \sin(\alpha - 2 \theta + \delta) & = & -\sin \delta + \frac{1}{4}\sin (2 \theta - \delta) \label{cons2}
\end{eqnarray}
\end{subequations}
\noindent Let
$$ \epsilon_1 \ \doteq \ \sin(\alpha - 2 \theta + \delta) - \sin(\alpha - 2 \theta), \quad \epsilon_2 \ \doteq \ \sin(2 \theta - \delta) - \sin( 2 \theta).$$
\noindent Using these definitions and substituting (\ref{cons1}) into (\ref{cons2}) we obtain
\begin{eqnarray}
\frac{1}{2.5} \left[ \sin \theta + \sin \delta + \epsilon_1 \right] & = & -\sin \delta + \frac{1}{4} \left[ \sin (2 \theta) + \epsilon_2 \right], \quad \mbox{or} \\
\Delta(\theta) & = & \frac{7}{2} \sin \delta + \epsilon_1 - \frac{5}{8} \epsilon_2. \label{cons3}
\end{eqnarray}
\noindent We observe that for any angle $\phi$ with 
$$ \max\{0, \pm \delta\}  \ \le \ \phi \ \le \ \min \{\pi/2, \pi/2 \pm \delta\}, $$
we have
\begin{eqnarray}
&& | \sin(\phi) - \sin(\phi - \delta) | \ \le \ |\delta|. \label{sincurve}
\end{eqnarray}
Then using (\ref{sincurve}) and (\ref{basic3}) (or (\ref{basic4})) we have $|\epsilon_1| \le \delta$ (resp., $|\epsilon_2| \le \delta$).  Hence, from (\ref{cons3}) and 
(\ref{basic3})
\begin{eqnarray}
| \Delta(\theta) | & \le & 0.02563.
\end{eqnarray}
\noindent {\bf 6.} As a consequence, using Lemma \ref{curvy} we have that either $\theta$ is close
to zero or close to $\theta_0$, or more precisely
\begin{eqnarray}
&& \mbox{either} \ \ \mbox{$0 \le \theta \le 0.1057$} \quad \mbox{or} \quad \mbox{$0.578 \le \theta \le 0.6952 \, (< \pi/4)$}.
\end{eqnarray}
We will refer to these two modes of operation as Mode I and II, respectively.\\

\noindent {\bf 7.} Define the throughput of the transmission system to be the
total flow sent out from bus $0$, which is the same as the total flow received
by bus $4$.  Hence the throughput equals
$$ \frac{1}{4} \sin(2\theta - \delta) + \sin(\theta).$$
As per observation {\bf 6} we now have that
\begin{itemize}
\item [(I)] In Mode I the throughput is less than $0.1592$.  
The solution using $\theta = 0.1057$ and $\delta = 0$ is feasible and attains
throughput greater than $0.1579$.
\item [(II)] In Mode II the throughput is at least $0.77464$ and less than $0.88671$.
The solution using $\theta = 0.6952$ and $\delta = 0$ is feasible and attains
throughput greater than $0.88648$.
\end{itemize}

\subsection{NP-hardness construction}
Here we consider the following problem:\\

\noindent THROUGHPUT: Given a transmission system in the lossless AC power 
flow model, with a single generator and a single load, and a value $T \ge 0$,
is there a feasible solution where at least $T$ units of power are transmitted from
the generator to the load?\\

\noindent We will show that this problem is strongly NP-hard using a reduction 
from one-in-three 3-SAT, defined as follows:\\

\noindent ONE-IN-THREE 3-SAT: Given clauses $C_1, \ldots, C_m$ on boolean 
variables $x_1, \ldots, x_n$, where each $C_i$ uses three literals, is there a truth assignment to the $x_j$
where in each clause there is precisely one true literal?\\

\noindent The construction will rely on several adaptations of the transmission system B considered in the previous section.  In what follows, we set
$$ S \, \doteq \, 0.1592, \quad H \, \doteq \, 0.8864.$$
\noindent We assume we are given an instance in one-in-three 3SAT with notation as above. For each variable $x_j$, we first construct a ``variable'' 
network $V(j)$ using two copies of the B network, denoted $B_j$ and $\bar B_j$,
respectively
plus three additional buses, $s_{j}, \, \bar s_{j}$ and $t_j$.  
\begin{wrapfigure}{r}{0.45\textwidth}
  \begin{center}
    \includegraphics[width=0.43\textwidth]{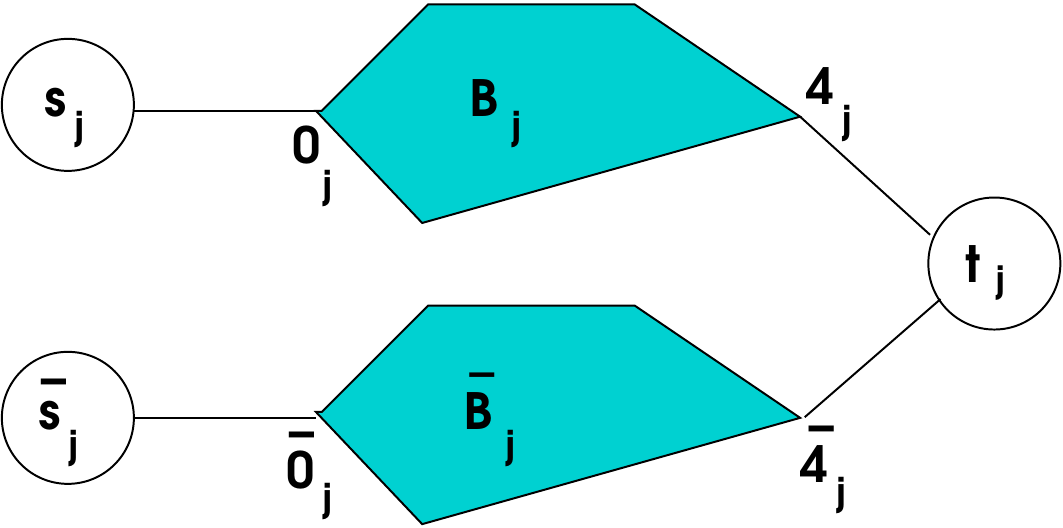}
  \end{center}
\vskip -15pt
  \caption{Network $V(j)$} \label{Vj}
\vskip 25pt
\end{wrapfigure}
\vskip 25pt
\noindent There is a line connecting bus $s_{j}$ (resp., $\bar s_{j}$) to the
copy of bus $0$ in $B_j$ ($\bar B_j$), denoted by $0_j$ ($\bar 0_j$).  Each of the two copies of bus 4 (denoted $4_j$ and $\bar 4_j$, respectively) is connected by a line to bus $t_j$.  See Figure \ref{Vj}.  The lines connecting 
$s_{j}, \, \bar s_{j}$ and $t_j$ to the $B$ networks have unit reactance and very large
limit.  

\begin{wrapfigure}[17]{r}{0.51\textwidth}
  \begin{center}
    \includegraphics[width=0.5\textwidth]{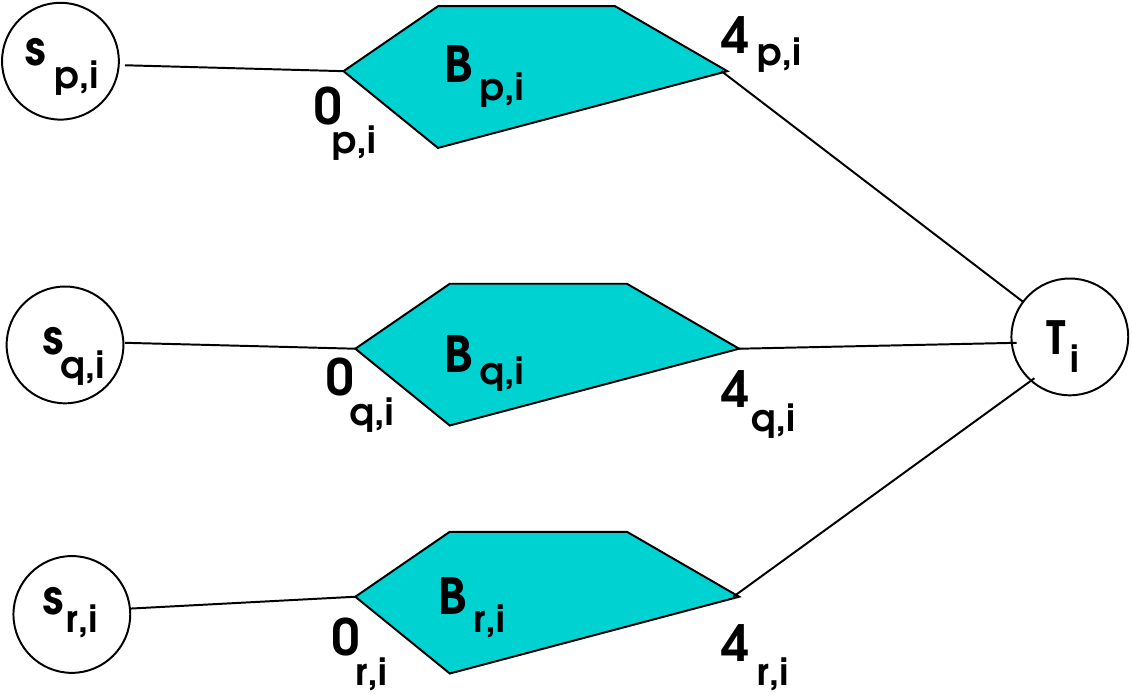}
  \end{center}
\vskip -15pt
  \caption{Network $C(i)$ corresponding to clause $C_i = (p \vee q \vee r)$.} \label{Ci}
\end{wrapfigure}

\noindent Next, for each clause $C_i$ we construct a ``clause'' network $C(i)$.
Suppose $C_i = (p \vee q \vee r)$.  The network $C(i)$ will contain three copies of the B network, denoted
$B_{p,i}$,  $B_{q,i}$ and $B_{r,i}$.  See Figure \ref{Ci}. \\

The copy of bus $0$ in $B_{p,i}$ is labeled
$0_{p, i}$ and likewise with the other copies of bus $0$ and the copies of bus $4$.
There are additional buses $s_{p,i}$, $s_{q,i}$ and $s_{r,i}$, connected to
$0_{p,i}$, $0_{q,i}$, and $0_{r,i}$, respectively, and a bus $T_i$ connected to
$4_{p,i}$, $4_{q,i}$, and $4_{r,i}$.  These new lines have unit reactance and 
large limit.  

\begin{figure}[h] 
\centering
\includegraphics[width=4in]{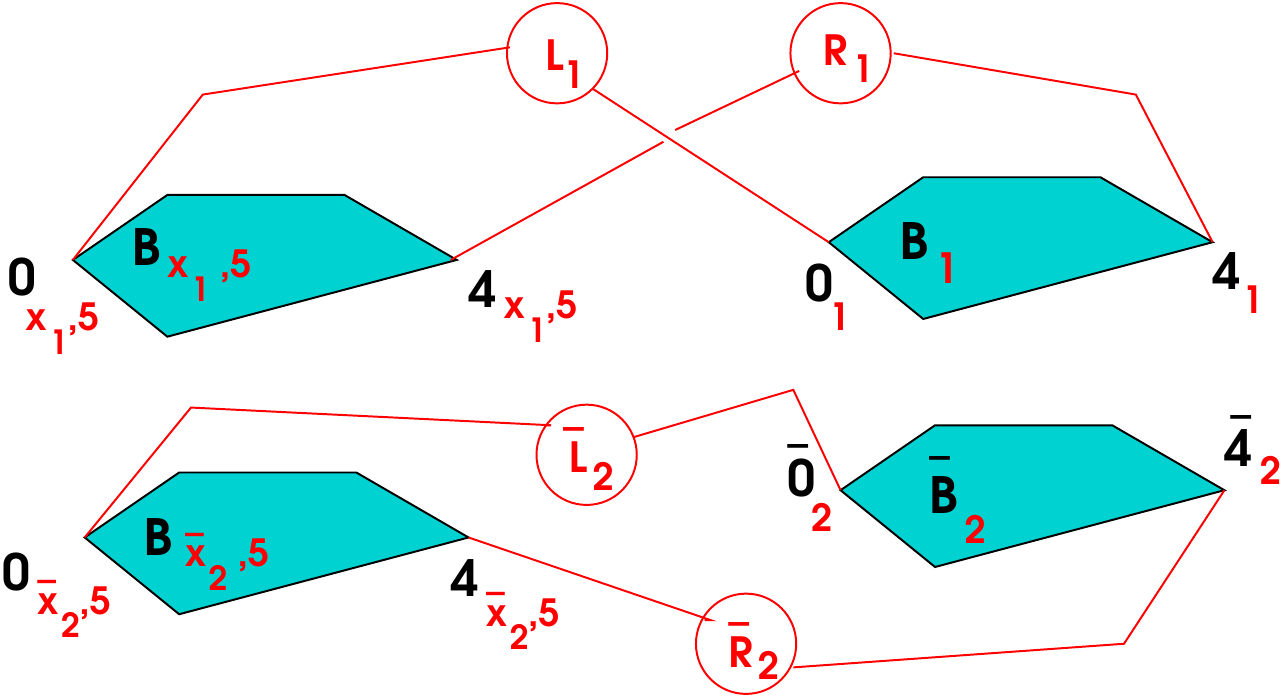}
\caption{Assembly in the case $C_5 \ = \ (\, x_1 \vee \bar x_2 \vee x_3 \,)$.}\label{assemble}
\end{figure}

\noindent The variable networks and clause networks are assembled together as follows.  
\begin{wrapfigure}[11]{r}{0.31\textwidth}
  \begin{center}
    \includegraphics[width=0.3\textwidth]{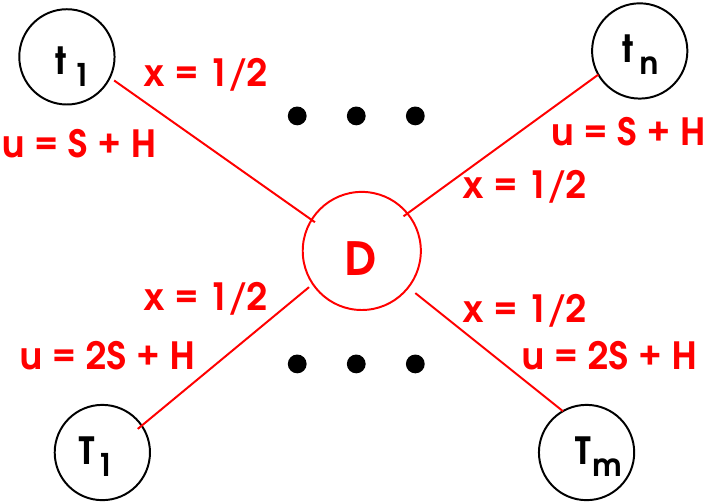}
  \end{center}
\vskip -15pt
  \caption{Attaching bus $D$} \label{D}
\end{wrapfigure}
\noindent (a) There is an additional bus, $D$. Bus $D$ is connected to all buses
$t_j$ ($1 \le j \le n$) with lines with reactance $1/2$ and limit $S + H$.
Bus $D$ is also connected to all buses
$T_i$ ($1 \le i \le m$) with lines with reactance $1/2$ and limit $2S + H$.  See 
Figure \ref{D}.\\
\noindent (b) For each variable $x_j$ ($1 \le j \le n$) we have four additional 
buses, denoted $L_j$, $R_j$, $\bar L_j$ and $\bar R_j$ respectively.  Bus
$L_j$ (resp., $R_j$) is connected to $0_j$ ($4_j$), and, for any clause $C_i$ containing $x_j$, bus $L_j$ (resp., $R_j$) is connected $0_{x_j, i}$ ($4_{x_j, i}$).
Likewise, bus
$\bar L_j$ (resp., $\bar R_j$) is connected to $\bar 0_j$ ($\bar 4_j$), and, for any clause $C_i$ containing $\bar x_j$, bus $\bar L_j$ (resp., $\bar R_j$) is connected $\bar 0_{x_j, i}$ ($\bar 4_{x_j, i}$).  All lines mentioned here have limit $1/20$ and unit reactance.  See Figure \ref{assemble} for an example.

\noindent (c) Bus $D$ is the only load.
For $1 \le j \le n$, buses $s_j$ and $\bar s_j$ are generators, each with capacity $S + H$.  $1 \le i \le m$, let $C_i = (p \vee q \vee r)$.  Then buses $s_{p,i}, \, s_{q,i}$ and $s_{r,i}$ are generators, each with capacity $2S + H$.  \\

The following two lemmas establish the NP-hardness result.
\begin{lemma} \label{first} Suppose there is a feasible solution where the total demand
consumed at node $D$ is at least $n(S + H) \, + \, m(2S + H)$. Then the 
instance of one-in-three 3SAT is satisfiable.
\end{lemma}
\noindent {\em Proof.} First we note that each of the subnetworks $B_j$, $\bar B_j$(for each variable $x_j$)  and $B_{w,i}$ (for each clause $C_i$) must operate
in Mode I or II as detailed above.  Further, by choice of the limits on the lines connected to bus $D$ it follows that each such line is operating at its limit.  
It follows that (a) for $1 \le j \le n$ one of $B_j$ and $\bar B_j$ operates in Mode I and the other in Mode
II, and (b) for $1 \le i \le m$ if $C_i = (p \vee q \vee r)$ then two of $B_{p,i}$, $B_{q,i}$ and $B_{r,i}$ operate in Mode I and the remaining one in Mode II.  

Further, suppose $x_j$ is one of the literals in clause $C_i$.  Since the limits
of lines $(0_j, L_j)$ and $(L_j, 0_{x_j,i})$ are $1/20$, and both lines have
reactance $1/2$, it follows that the absolute value of the difference of phase
angles at $0_j$ and $0_{x_j,i}$ is at most $2 \sin^{-1}(1/10) < 0.201$.  
Likewise, the absolute value of the difference of phase
angles at $4_j$ and $4_{x_j,i}$ is at less than $0.201$.  Hence $B_j$ operates 
in Mode I (or Mode II) if and only if $B_{x_j, i}$ operates in Mode I (Mode II, 
respectively).  

Similarly, if $\bar x_j$ is one of the literals in clause $C_i$, then 
$\bar B_j$ operates 
in Mode I (or Mode II) if and only if $B_{x_j, i}$ operates in Mode I (Mode II, 
respectively).  

The result is now established by using the truth assignment $x_j = true$ iff
$B_j$ operates in Mode II. \QED
\begin{lemma} Suppose the 
instance of one-in-three 3SAT is satisfiable. Then there is a feasible solution where the total demand
consumed at node $D$ is at least $n(S + H) \, + \, m(2S + H)$. 
\end{lemma}
\noindent {\em Proof.} Similar to that of Lemma \ref{first}. For $1 \le j \le n$
we operate $B_j$ operate in Mode II iff $x_j = true$, and for each clause 
$C_i = (p \vee q \vee r)$, if, say $p = x_j$ then we operate $B_{x_j, i}$ in the
same mode as $B_j$ and with buses $L_j$, $0_j$ and $0_{x_j, i}$ at the same phase
angle and buses $R_j$, $4_j$ and $4_{x_j, i}$ at the same phase
angle. And if $p = \bar x_j$ then we operate $B_{\bar x_j, i}$ in the
same mode as $\bar B_j$ and with a corresponding setting of phase angles.  This
setting of phase angles yields a feasible solution with the desired throughput,
where all lines incident with buses $L_j$ and $R_j$ carry zero flow. \QED
\subsection{Membership in NP}
The above proof shows that testing feasibility for a system of type (\ref{lossless}) is an NP-hard problem.  To prove NP-completeness we would need to argue for 
membership in NP.  A straightforward proof of such a fact, if true, is unlikely, for the reason that in a feasible solution very likely the $f_{ij}$ (and possibly
even some of the $\theta_i$) would be irrational values.

We conjecture that an approximate version of system (\ref{lossless}) where 
equation (\ref{ohm}) is replaced with $|\sin(\theta_i - \theta_j) \ - \ x_{ij} f_{ij}| \le \epsilon$ (where $\epsilon > 0$ is part of the input) belongs to NP.
\bibliographystyle{siam}
\bibliography{bananas}
\hspace{.1in}\\
\hspace{.1in}\\
\tiny Tue.Apr..9.125829.2019@littleboy

\end{document}